\documentstyle[12pt]{article}
\def\stern{\mbox{\LARGE $^*$}}
\begin{document}
\includegraphics{zibheader.eps}
\vspace*{6.5cm}
\begin{center}
{\Large Wolfram Koepf}
\vspace*{1cm}\\
{\Large Dieter Schmersau\stern}
\vspace*{2cm}\\
\Large{{\bf  On the De Branges Theorem}}
\end{center}
\vfill
\vfill
\vfill
$*$ Fachbereich Mathematik der Freien Universit\"at Berlin
\vfill
\hrule
\vspace*{3mm}
Preprint SC 95--10 (March 1995)

\thispagestyle{empty}
\setcounter{page}{0}
\eject
\hoffset -1cm
\footskip=1cm
\parindent=0pt
\title{On the De Branges Theorem}
\date{}
\author{Wolfram Koepf\\
Dieter Schmersau\\
email: {\tt koepf@zib-berlin.de}
}
\maketitle
%
%
%
%
\newcommand{\hypergeom}[5]{\mbox{$
_#1 F_#2\left.
\!\!
\left(
\!\!\!\!
\begin{array}{c}
\multicolumn{1}{c}{\begin{array}{c}
#3
\end{array}}\\[1mm]
\multicolumn{1}{c}{\begin{array}{c}
#4
            \end{array}}\end{array}
\!\!\!\!
\right| \displaystyle{#5}\right)
$}
}
\setlength{\textheight}{8.5in}
\setlength{\textwidth}{6.5in}
\setlength{\evensidemargin} {0in}
\setlength{\oddsidemargin} {0in}
\setlength{\topmargin} {0in}

\newcommand{\punktpunkt}[1]{\stackrel{\mbox{..}}{#1}\!}

\newcommand{\RR}{{\rm I\! R}}
\newcommand{\DD}{{\rm I\! D}} 
\newcommand{\CC}{{\; \rm l\!\!\!	C}}
\newcommand{\NN}{{\rm I\!I\!\! N}}
\newcommand{\ZZ}{\Bbb Z}
\newcommand{\C}{{\rm {\mbox{C{\llap{{\vrule height1.5ex}\kern.4em}}}}}}
\newcommand{\N} {{\rm {\mbox{\protect\makebox[.15em][l]{I}N}}}}
\renewcommand{\H} {{\rm {\mbox{\protect\makebox[.15em][l]{I}H}}}}
\newcommand{\Q} {{\rm {\mbox{Q{\llap{{\vrule height1.5ex}\kern.5em}}}}}}
\newcommand{\R} {{\rm {\mbox{\protect\makebox[.15em][l]{I}R}}}}
\newcommand{\D} {{\rm {\mbox{\protect\makebox[.15em][l]{I}D}}}}
\newcommand{\Z} {{\rm {\mbox{\protect\makebox[.2em][l]{\sf Z}\sf Z}}}}
\newcommand{\Rp}
{{\rm {\mbox{$\mbox{\protect\makebox[.15em][l]{I}R}^{\scriptscriptstyle+}$}\index{R+@\mbox{$\mbox{\protect\makebox[.15em][l]{I}R}^{\scriptscriptstyle+}$}, 
Notation for the positive real numbers}}}}
\newcommand{\Rm}
{{\rm {\mbox{$\mbox{\protect\makebox[.15em][l]{I}R}^{\scriptscriptstyle-}$}\index{R-@\mbox{$\mbox{\protect\makebox[.15em][l]{I}R}^{\scriptscriptstyle-}$}, 
Notation for the negative real numbers}}}}
\newcommand{\Pt}{\widetilde{P}}
\newcommand{\Bt}{\widetilde{B}}
\newcommand{\Cd}{\widehat{{\rm \;l\!\!\! C}}}
\newcommand{\bq}[1]{|#1|^2}
\newcommand{\vf}[1]{(1-\bq{#1})}
\newcommand{\vfq}[1]{\vf{#1}^2}
\newcommand{\Bf}[1]{\vf{#1}\left|\frac{f''}{f'}(#1)\right|}
\newcommand{\Nf}[1]{\vfq{#1}|S_f(#1)|}
\newcommand{\Kf}[1]{\left|-\kon{#1}+\frac{1}{2}\vf{#1}\frac{f''}{f'}(#1)\right|}
\newcommand{\ed}[1]{\frac{1}{#1}}
\newcommand{\aut}[1]{\frac{{\textstyle{z+#1}}}{{\textstyle{1+\kon{#1}z}}}}
\newcommand{\au}[2]{#1\aut{#2}}
\newcommand{\subs}[2]{\left. \makebox{\rule{0in}{2.5ex}} #1 \rb_{#2}}
\newcommand{\subst}[3]{\left. \makebox{\rule{0in}{2.5ex}} #1 \rb_{#2}^{#3}}
\newcommand{\funkdef}[3]{\left\{\begin{array}{ccc}
                                #1 && \mbox{\rm{if} $#2$ } \\
                                #3 && \mbox{\rm{otherwise}}
                                \end{array}  
                         \right.}       
\newcommand{\funkdeff}[4]{\left\{\begin{array}{ccc}
                                 #1 && \mbox{\rm{if} $#2$ } \\
                                 #3 && \mbox{\rm{if} $#4$ } 
                                 \end{array}
                          \right.}
\newcommand{\funkdefff}[6]{\left\{\begin{array}{ccc}
                                 #1 && \mbox{{if} $#2$ } \\
                                 #3 && \mbox{{if} $#4$ } \\
                                 #5 && \mbox{{if} $#6$ }
                                 \end{array}
                          \right.}
\newcommand{\ueber}[2]{
                       \Big( \!
                       {{\small
                       \begin{array}{c}
                          #1\\
                          #2
                          \end{array}
                       }}
                       \! \Big) }
\newcommand{\function}[4]{
                          \begin{array}{rcl}#1&\pf&#2\\
                          #3&\mapsto &#4
                          \end{array} }
\newcommand{\pr}{\vspace{-2mm}\absatz{{\sl Proof:}}\hspace{5mm}}
\newcommand{\eop}{\hfill$\Box$\par\medskip\noindent}
\newcommand{\absatz}{\par\medskip\noindent}
\renewcommand{\Re}{{\rm Re\:}}
\renewcommand{\Im}{{\rm Im\:}}
\newcommand{\co}{{\rm co\:}}
\newcommand{\coq}{\overline{{\rm co}} \:}
\newcommand{\ex}{{\rm E\:}}
\newcommand{\In}{\in}
\newcommand{\ro}{\varrho}
\newcommand{\om}{\omega}
\newcommand{\al}{\alpha}
\newcommand{\bb}{\beta}
\newcommand{\la}{\lambda}
\newcommand{\eps}{\varepsilon}
\newcommand{\ph}{\varphi}
\renewcommand{\phi}{\varphi}
\newcommand{\si}{\sigma}
\newcommand{\ka}{\varkappa}
\newcommand{\th}{\theta}
\newcommand{\g}{\gamma}
\newcommand{\de}{\partial}
\newcommand{\fD}{f(\D)}
\newcommand{\sumi}{\sum\limits_{k=0}^{\infty}}
\newcommand{\sumei}{\sum\limits_{k=1}^{\infty}}
\newcommand{\union}{\bigcup\limits_{k=1}^{n}}
\newcommand{\sumn}{\sum\limits_{k=1}^{n}}
\newcommand{\prodn}{\prod\limits_{k=1}^{n}}
\newcommand{\intd}{\int\limits_{\de\DD}}
\newcommand{\menge}[3]{\left\{#1 \In #2 \; \lb \; #3 \right. \right\} }
\newcommand{\mk}{\mu_{k}}
\newcommand{\xk}{x_k}
\newcommand{\yk}{y_k}
\newcommand{\xn}{x_n}
\newcommand{\yn}{y_n}
\newcommand{\ak}{\al_k}
\newcommand{\bk}{\bb_k}
\newcommand{\kn}{\mbox{$(k=1, \ldots ,n)$}}
\newcommand{\kno}{\mbox{$k=1, \ldots ,n$}}
\newcommand{\sub}{\prec}
\newcommand{\ld}[1]{\frac{f''}{f'}(#1)} 
\newcommand{\limn}{\lim\limits_{n\rightarrow\infty}}
\newcommand{\lsz}{\limsup\limits_{z\rightarrow\de\DD}}
\newcommand{\limr}{\lim\limits_{r\rightarrow 1}}
\newcommand{\liz}{\liminf\limits_{z\rightarrow\de\DD}}
\newcommand{\supD}[1]{\sup\limits_{#1\In \DD}}
\newcommand{\infD}[1]{\inf\limits_{#1\In \DD}}
\newcommand{\maxn}{\max\limits_{1 \leqq k \leqq n}}
\newcommand{\minn}{\min\limits_{1 \leqq k \leqq n}}
\newcommand{\ord}{{\rm ord\:}}
\newcommand{\gleich}[1]{\stackrel{{{\rm (#1)}}}{\longeq}}
\newcommand{\folgt}[1]{\stackrel{{{\rm (#1)}}}{\Pf}}
\newcommand{\gegen}[1]{\stackrel{{{\rm (#1)}}}{\pf}}
\newcommand{\nach}[2]{(#1)$\dpf$(#2):}
\newcommand{\Exp}{\subset\!\subset}
\newcommand{\lleq}{\stackrel{_{{\scriptscriptstyle \Exp}}}
           {_{{\scriptscriptstyle \sim}}}}
\newcommand{\Sub}{{\rm Sub\:}}
\newcommand{\norm}[2]{\frac{#1\circ #2-#1\circ #2(0)}{(#1\circ #2)'(0)}} 
\renewcommand{\dim}{{\rm{dim}}_{_{H^{p}}}}
\renewcommand{\ll}{<\!<}
\newcommand{\1}{{\bf{1}}}
\newcommand{\2}{{\bf{2}}}
\newcommand{\3}{{\bf{3}}}
\newcommand{\4}{{\bf{4}}}
\newcommand{\5}{{\bf{5}}}
\newcommand{\6}{{\bf{6}}}
\newcommand{\7}{{\bf{7}}}
\newcommand{\8}{{\bf{8}}}
\newcommand{\9}{{\bf{9}}}
\newcommand{\0}{{\bf{0}}}
\def\finis{\hbox{$\bigtriangleup$}}    


\newcommand{\luc}{locally uniform convergence}
\newcommand{\hp}{half\-plane}
\newcommand{\sq}{sequence}
\newcommand{\an}{analytic}
\newcommand{\af}{analytic\ function}
\newcommand{\Ne}{Ne\-ha\-ri\ ex\-pression}
\newcommand{\Ke}{Koe\-be\ ex\-pression}
\newcommand{\Be}{Becker\ ex\-pression}
\newcommand{\fc}{function}
\newcommand{\uv}{univalent}
\newcommand{\uf}{univalent\ function}
\newcommand{\SC}{Schwarz-Chri\-stof\-fel}
\newcommand{\pol}{polyg\-on}
\newcommand{\cv}{con\-vex}
\newcommand{\ctc}{close-to-con\-vex}
\newcommand{\Ck}{Ca\-ra\-th\'eo\-dory\ kernel}
\newcommand{\st}{sector}
\newcommand{\sm}{similar}
\newcommand{\br}{boundary\ rotation}
\newcommand{\bbr}{bounded\ \br}
\newcommand{\KM}{Krein-Mil\-man}


\newcommand{\til}{\widetilde}
\newcommand{\pf}{\rightarrow}
\newcommand{\Pf}{\;\;\;\longrightarrow\;\;\;}
\newcommand{\dpf}{\Rightarrow}
\newcommand{\Dpf}{\;\;\;\Longrightarrow\;\;\;}
\newcommand{\kon}{\overline}
\newcommand{\be}{\begin{equation}}
\newcommand{\ee}{\end{equation}}
\newcommand{\bea}{\begin{eqnarray}}
\newcommand{\eea}{\end{eqnarray}}
\newcommand{\beao}{\begin{eqnarray*}}
\newcommand{\eeao}{\end{eqnarray*}}
\newcommand{\lequiv}{\;\;\;\Longleftrightarrow\;\;\;}
\newcommand{\longeq}{\;\;\;=\!\!=\!\!=\;\;\;}
\newcommand{\leqq}{\leq}
\newcommand{\geqq}{\geq}
\newcommand{\gl}{\;\leftrightarrow\;}
\newcommand{\lk}{\left(}
\newcommand{\rk}{\right)}
\newcommand{\lb}{\left|}
\newcommand{\rb}{\right|}
\newcommand{\bi}{\bibitem}


\newcommand{\bT}{\begin{theorem}}
\newcommand{\eT}{\end{theorem}}
\newcommand{\bL}{\begin{lemma}}
\newcommand{\eL}{\end{lemma}}
\newcommand{\bC}{\begin{corollary}}  
\newcommand{\eC}{\end{corollary}}
\newcommand{\bt}{\begin{tabbing} 12345 \= \kill}
\newcommand{\et}{\end{tabbing}}


\hyphenation{qua-si-disk qua-si-circle 
non-smooth
pa-ram-e-trized pa-ram-e-tri-zation
geo-met-ric
}


\newcommand{\bbegin}{{\bf{begin}}}
\newcommand{\eend}{{\bf{end}}}
\newcommand{\iif}{{\bf{if}}}
\newcommand{\tthen}{{\bf{then}}}
\newcommand{\wwhile}{{\bf{while}}}
\newcommand{\ddo}{{\bf{do}}}
\newcommand{\ffor}{{\bf{for}}}
\newcommand{\sstep}{{\bf{step}}}
\newcommand{\llet}{{\bf{let}}}
\newcommand{\pprocedure}{{\bf{procedure}}}
\newcommand{\aand}{{\bf{and}}}
\newcommand{\nnot}{{\bf{not}}}
\newcommand{\oor}{{\bf{or}}}
\newcommand{\lllet}{{\bf{let}}}
\newcommand{\eexit}{{\bf{exit}}}
\newcommand{\rreturn}{{\bf{return}}}
\newcommand{\uuntil}{{\bf{until}}}

\newcommand{\abs}{\\[3mm]}

\newtheorem{theorem}{Theorem}
\newtheorem{algorithm}{Algorithm}
\newtheorem{lemma}{Lemma}
\newtheorem{corollary}{Corollary}
\newtheorem{definition}{Definition}
\begin{abstract}
Recently, Todorov and Wilf independently 
realized that de Branges' original proof of the
Bieberbach and Milin conjectures and the proof that was later given by
Weinstein deal with the same special function system that de Branges had
introduced in his work.

In this article, we present an elementary proof of this statement
based on the defining differential equations system rather than the
closed representation of de Branges' function system. Our proof does neither use
special functions (like Wilf's) nor the residue theorem (like Todorov's)
nor the closed representation (like both), but is purely algebraic. 

On the other hand, by a similar algebraic treatment, the closed representation
of de Branges' function system is derived.
In a final section, we give a simple representation of a generating 
function of the de Branges functions.

Our whole contribution can be looked at as the study of properties of
the Koebe function.
Therefore, in a very elementary manner it is shown that
the known proofs of
the Bieberbach and Milin conjectures can be understood as a consequence
of the L\"owner differential equation, plus properties of the Koebe function.
\end{abstract}
 
\section{Introduction}

Let $S$ denote the family of analytic and univalent functions
$f(z)=z+a_2 z^2+\ldots$ of the unit disk $\D$. $S$ is compact with
respect to the topology of locally uniform convergence so that
\mbox{$k_{n}:=\max\limits_{{f\in S}}|a_{n}(f)|$} exists. In 1916
Bieberbach \cite{Bieberbach}
proved that $k_2=2$, with equality if and only if $f$ is a rotation of
the {\sl Koebe function}
\be
K(z):=\frac{z}{(1-z)^{2}}=\ed{4}\lk\lk\frac{1+z}{1-z}\rk^{2}-1\rk=
\sum\limits_{n=1}^{\infty}nz^{n}\;,
\label{eq:Koebe function}
\ee
and in a footnote he mentioned ``Vielleicht ist
\"uberhaupt $k_{n}=n$.''. This statement is known as the {\sl Bieberbach
conjecture.}

In 1923 L\"owner \cite{Loewner2} proved the Bieberbach conjecture for
$n=3$. His method was to embed a univalent function $f(z)$ into a {\sl L\"owner
chain\/}, i.e.\ a family $\left\{f(z,t) \; \lb \; t\geqq 0 \right.
\right\}$ of univalent functions of the form
\[
f(z,t)=e^{t}z+\sum\limits_{n=2}^{\infty}a_{n}(t)z^{n},
\;\;\;(z\in\D, t\geqq 0, a_{n}(t)\in\C\;(n\geqq 2))
\]
which start with $f$
\[
f(z,0)=f(z)\;,
\]
and for which the relation
\be
\Re p(z,t)=\Re\lk\frac{{\dot{f}}(z,t)}{z f'(z,t)}\rk>0
\;\;\;\;\;\;\;\;\;(z\in\D)
\label{eq:p(z,t)}
\ee
is satisfied. Here $'$ and $\dot{}\:$ denote the partial derivatives
with respect to $z$ and  $t$, respectively. Equation~(\ref{eq:p(z,t)})
is referred to as the
{\sl L\"owner differential equation\/}, and
geometrically it states that the image domains of
$f_t$ expand as $t$ increases.

The history of the Bieberbach conjecture showed that it was easier to
obtain results about the {\sl logarithmic coefficients\/} of a univalent
function $f$, i.e.\ the coefficients $d_n$ of the expansion
\[
\phi(z)=\ln\frac{f(z)}{z}=:\sum\limits_{n=1}^{\infty}d_{n}z^{n}
\]
rather than for the coefficients $a_n$ of $f$ itself. So Lebedev and
Milin \cite{LM}
in the mid sixties developed methods to exponentiate such
information. They proved that if
for $f\in S$ the {\sl Milin conjecture\/}
\[
\sum\limits_{k=1}^{n}(n+1-k)\lk k|d_{k}|^{2}-\frac{4}{k}\rk\leqq 0
\]
on its logarithmic coefficients is satisfied for some $n\in\N$, then 
the Bieberbach conjecture for the index $n+1$ follows.

In 1984 de Branges \cite{Bra}
verified the Milin, and therefore the Bieberbach
conjecture, and in 1990, Weinstein \cite{Weinstein} gave a different proof.
A reference concerning de Branges' proof is \cite{FP},
and a German language summary of the history of the Bieberbach conjecture
and its proofs was given in \cite{Koepf}.

Both proofs use special function systems, and
independently, Todorov \cite{Todorov} and Wilf \cite{Wilf}
showed that these essentially are the same.

In this article, we present an elementary proof of this statement.
Our considerations are
based on the defining differential equations system rather than the
closed representation of de Branges' function system. Our proof neither uses
special functions (like Wilf's) nor the residue theorem (like Todorov's)
nor the closed representation (like both), but is purely algebraic.

On the other hand, by a similar algebraic treatment, the closed representation
of de Branges' function system is derived.

In a final section, we give a simple representation of a generating 
function of the de Branges functions.

Our whole contribution can be looked at as the study of properties of
the Koebe function.
Therefore, in a very elementary manner it is shown 
that the known proofs of
the Bieberbach and Milin conjectures can be understood as a consequence
of the L\"owner differential equation, plus properties of the Koebe function.

\section{The L\"owner Chain of the Koebe Function}
\label{sec:The Lowner Chain of the Koebe Function}

In this section, we consider the L\"owner chain
\be
w(z,t):=K^{-1}\Big(e^{-t}K(z)\Big)\quad(z\in\D, t\geq 0)
\label{eq:Koebe}
\ee
of bounded univalent functions in the unit disk $\D$
which is defined in terms of the Koebe function (\ref{eq:Koebe function}).
Since $K$ maps the unit disk onto the entire plane slit along the
negative $x$-axis in the interval $(-\infty,1/4]$,
$w(\D,t)$ is the unit disk with a radial slit increasing with $t$.
The function $w(z,t)$ is implicitly given by the equation
\be
K(z)=e^t\,K(w(z,t))
\label{eq:w(z,t)-implicitly}
\;,
\ee
and satisfies the L\"owner type differential equation (we omit the arguments)
\be
\dot w=-\frac{1-w}{1+w}\,w
\label{eq:w(z,t)LoewnerDE}
\ee
(compare e.\ g.\ \cite{Pommerenke}, Chapter 6) which is obtained 
differentiating (\ref{eq:w(z,t)-implicitly}) with respect to $t$
\be
0=e^t\,K(w)+e^t\,K'(w)\,\dot w
\;,
\label{eq:diffwrtt}
\ee
hence
\[
\dot w=-\frac{K(w)}{K'(w)}=-\frac{w}{(1-w)^2}\frac{(1-w)^3}{1+w}=
-\frac{1-w}{1+w}\,w
\;.
\]
In this section, we deduce a closed representation of the Taylor
coefficients $A_n(t)$ of
\be
w(z,t)=\sum_{n=1}^\infty A_n(t)\,z^n
\;.
\label{eq:A(k,t)}
\ee
In particular, by the normalization of the Koebe function, we have
$\subs{\frac{K(z)}{z}}{z=0}=1$, hence by (\ref{eq:w(z,t)-implicitly})
\[
\frac{K(z)}{z}=e^t\,\frac{K(w(z,t))}{z}=
e^t\,\frac{K(w(z,t))}{w(z,t)}\frac{w(z,t)}{z}
\]
and letting $z\pf 0$ therefore gives $A_1(t)=e^{-t}$.

To deduce the general result in an elementary way
(for a shorter deduction using Gegenbauer polynomials, see 
\S~\ref{sec:Connection with the Gegenbauer Polynomials}),
we begin with some
lemmas. The first lemma states a linear partial differential equation 
(different from the nonlinear L\"owner differential equation
(\ref{eq:w(z,t)LoewnerDE})) for $w(z,t)$:
\bL[Partial differential equation]
{\rm
The function $w(z,t)$ satisfies the linear partial differential equation
\be
(z-1)zw'(z,t)=(z+1)\dot w(z,t)
\label{eq:partial differential equation}
\ee
with the initial function
\[
w(z,0)=z
\;.
\]
}
\eL
\pr
Differentiating (\ref{eq:w(z,t)-implicitly}) with respect to both $z$ and $t$
yields the equations
\[
K'(z)=e^t\,K'(w(z,t))\,w'(z,t)
\]
and (\ref{eq:diffwrtt}), from which we deduce
\[
\frac{zw'(z,t)}{\dot w(z,t)}=
-\frac{zK'(z)}{e^t\,K'(w(z,t))}
\frac{K'(w(z,t))}{K(w(z,t))}
=
-\frac{zK'(z)}{e^t\,K(w(z,t))}=
-\frac{zK'(z)}{K(z)}
=
-\frac{1+z}{1-z}
\]
where we used (\ref{eq:w(z,t)-implicitly}) once again,
and (\ref{eq:Koebe function}).
The initial function is determined trivially.
\eop
As a consequence we have for the coefficients $A_n(t)$ of
$w(z,t)$:
\bL
{\bf (Differential equations system for coefficient functions)}
\label{l:Differential equations system for coefficient functions}
{\rm
The coefficients $A_n(t)$ satisfy the system of linear differential equations
\be
(n-1)\,A_{n-1}(t)-n\,A_n(t)=
\dot A_{n-1}(t)+\dot A_{n}(t)
\;,
\quad\quad A_n(0)=0
\quad\quad (n\geq 2)
\label{eq:ARekursion}
\ee
and
\be
-A_1(t)=\dot A_1(t)\;,
\quad\quad A_1(0)=1
\;.
\label{AAnfangsrekursion}
\ee
}
\eL
\pr
This follows directly by summing (\ref{eq:partial differential equation})
for $n=0,\ldots,\infty$, and equating coefficients.
\eop
Starting with the solution $A_1(t)=e^{-t}$ of (\ref{AAnfangsrekursion}),
by induction we see that $A_n(t)$ is a polynomial of degree $n$ in $e^{-t}$.
Therefore we may introduce the variable $y:=e^{-t}$, and define the
polynomials $B_n(y)$ by
\be
A_n(t)=B_n(y)=B_n(e^{-t})
=\sum_{j=1}^n a_j^{(n)}\,e^{-jt}
=\sum_{j=1}^n a_j^{(n)}\,y^j
\;,
\label{eq:ajk}
\ee
so that in terms of $B_n(y)$,
Lemma~\ref{l:Differential equations system for coefficient functions}
reads as follows:
\bL
{\bf (Differential equations system for coefficient functions)}
{\rm
The functions $B_n(y)$ satisfy the system of linear differential equations
\be
y\,(B'_{n}(y)+B'_{n-1}(y))
=
n\,B_n(y)-(n-1)\,B_{n-1}(y)
\;,
\quad\quad B_n(1)=0
\quad\quad (n\geq 2)
\label{eq:BRekursion}
\;.
\ee
}
\eL
For the numbers $a_j^{(n)}$, we deduce
\bL[Recurrence equations]
{\rm
For the numbers $a_j^{(n)}$ defined by (\ref{eq:ajk}), the simple recurrence
equation
\be
(n-j)a_j^{(n)}=(n-1+j)a_j^{(n-1)}
\quad\quad (1\leq j \leq n-1, n\geq 2)
\label{eq:RE1 ajk}
\ee
is valid. Therefore, we have
\be
a_j^{(n)}=\ueber{n+j-1}{n-j}a_j^{(j)}
\quad\quad (1\leq j \leq n, n\geq 2)
\label{eq:RE3 ajk}
\ee
and the initial value
\[
a_1^{(n)}=n
\;.
\]
}
\eL
\pr
For $j=n$, Equation (\ref{eq:RE3 ajk}) is trivial. Therefore assume
$1\leq j \leq n-1, n\geq 2.$
Substituting (\ref{eq:ajk}) into (\ref{eq:BRekursion}),
and equating coefficients of $y^j\;
(1\leq j \leq n-1)$ results in (\ref{eq:RE1 ajk}). 
 From (\ref{eq:RE1 ajk}), we get the telescoping product
\[
a_j^{(n)}=\frac{n-1+j}{n-j}a_j^{(n-1)}=
\frac{(n-1+j)(n-2+j)\cdots (2j)}{(n-j)(n-j-1)\cdots 1}a_j^{(j)}=
\ueber{n-1+j}{n-j}a_j^{(j)}
\;,
\]
and hence (\ref{eq:RE3 ajk}).

Using $A_1(t)=e^{-t}$, we get $a_1^{(1)}=1$, so that by
(\ref{eq:RE3 ajk}), we finally have
\[
a_1^{(n)}=\ueber{n}{n-1}\,a_1^{(1)}=n
\;,
\]
and we are done.
\eop
Our next step is to derive an ordinary differential equation valid for
$B_n(y)$:
\bL
{\bf (Ordinary differential equation for coefficient functions)}
\label{l:Ordinary differential equation for coefficient functions}
{\rm
The function $B_n(y)\;(n\geq 1)$ satisfies the ordinary differential equation
\be
y^2\,(1-y)\,B_n''(y)+y\,(1-y)\,B_n'(y)+(n^2\,y-1)\,B_n(y)=0
\label{eq:ODE B}
\;.
\ee
}
\eL
\pr
For $n=1$, the statement is true, so assume $n\geq 2$.
We consider the function
\be
\Delta_n(y):=
y^2\,(1-y)\,B_n''(y)+y\,(1-y)\,B_n'(y)+(n^2\,y-1)\,B_n(y)
\;,
\label{eq:Deltakdefinition}
\ee
and show in a first step that $\Delta_n(y)$ satisfies the relation
\be
y\lk \Delta_n'(y)+\Delta_{n-1}'(y)\rk-n\Delta_n(y)+(n-1)\Delta_{n-1}(y)
=0
\;.
\label{eq:ODE DE}
\ee
(In other words, we show that $\Delta_n(y)$ satisfies the same system
of differential equations 
(\ref{eq:BRekursion}) as $B_n(y)$.)

To prove (\ref{eq:ODE DE}), 
we first solve (\ref{eq:BRekursion}) for $B_{n-1}'(y)$:
\[
B_{n-1}'(y)=\frac{n}{y}B_{n}(y)-\frac{n-1}{y}B_{n-1}(y)-B_{n}'(y)
\;.
\]
We take this equation and the first two derivatives thereof as replacement
rules for any occurrence of $B_{n-1}'(y)$, $B_{n-1}''(y)$, and $B_{n-1}'''(y)$
in the left hand side of (\ref{eq:ODE DE}). The resulting term reduces
to zero. This procedure can be easily done with the aid of a computer 
algebra system, and we leave these elementary algebraic transformations 
to the reader.

Therefore, by (\ref{eq:ODE B})--(\ref{eq:Deltakdefinition}) we have
$\Delta_1(y)\equiv 0$, and further $\Delta_n(1)=0$ for all $n\in\N$.
 From the induction hypothesis $\Delta_{n-1}(y)\equiv 0$, we get
the initial value problem
\[
y\Delta_n'(y)-n\Delta_n(y)=0
\;,
\quad\quad\mbox{and}\quad\quad
\Delta_n(1)=0
\;,
\]
and by integration the unique solution $\Delta_n(y)\equiv 0$ is deduced.
\eop
Obviously there is a corresponding ordinary differential equation for
$A_n(t)$, namely
\[
\left( 1 - {e^t} \right) \,\punktpunkt{A_n}(t)+
\left( {e^t} - {n^2} \right) \,A_n(t) 
=0
\;,
\]
which is simpler than (\ref{eq:ODE B}) in the sense that it does
not contain the first derivative explicitly, but which does not have
polynomial coefficients since $e^{t}$-terms occur.

As a consequence of the preceding lemmas, 
we find the following closed form representation of 
$a_j^{(n)}$:
\bT
{\bf (Coefficient representation of L\"owner chain of Koebe function)}
\label{th:Coefficient representation of Lowner chain of Koebe function}
{\rm
For the numbers $a_j^{(n)}$ defined by (\ref{eq:ajk}), we have the
closed form representation
\be
a_j^{(n)}=2(-1)^{j+1}\ueber{n+j-1}{n-j}\frac{(2j-1)!}{(j-1)!(j+1)!}
\quad\quad (1\leq j \leq n, n\geq 1)
\;.
\label{eq:ajk-formula}
\ee
Therefore, by (\ref{eq:ajk}), we have further
\be
A_n(t)=\sum_{j=1}^n 2(-1)^{j+1}\ueber{n+j-1}{n-j}\frac{(2j-1)!}{(j-1)!(j+1)!}
\,e^{-jt}
\quad\quad (n\geq 1)
\;,
\label{eq:Akfineal}
\ee
and finally by (\ref{eq:A(k,t)})
\be
w(z,t)=\sum_{n=1}^\infty
\sum_{j=1}^n 2(-1)^{j+1}\ueber{n+j-1}{n-j}\frac{(2j-1)!}{(j-1)!(j+1)!}
\,e^{-jt}
\,z^n
\;.
\label{eq:closedform:w(z,t)}
\ee
}
\eT
\pr
By (\ref{eq:RE3 ajk}), it remains to prove that
\be
a_j^{(j)}=
2(-1)^{j+1}\frac{(2j-1)!}{(j-1)!(j+1)!}
\;.
\label{eq:ajj}
\ee
Substituting (\ref{eq:ajk}) into 
the ordinary differential equation (\ref{eq:ODE B}), and equating coefficients
gives
\[
(n+1-j) (n-1+j) a_{j-1}^{(n)} + (j-1) (j+1) a_{j}^{(n)} = 0
\;,
\]
so that in particular for $n=j$
\[
a_{j}^{(j)} = -\frac{2j-1}{(j-1) (j+1)}\,a_{j-1}^{(j)}
=
-2\,\frac{2j-1}{j+1}\,a_{j-1}^{(j-1)}
\quad\quad (1\leq j \leq n-1, n\geq 2)
\;,
\]
by an application of (\ref{eq:RE1 ajk}), and therefore (\ref{eq:ajj}) follows
from $a_2^{(2)}=-2$. It is easily checked that (\ref{eq:ajk-formula})
remains true for $n=1$.
\eop

\section{Connection with the Gegenbauer Polynomials}
\label{sec:Connection with the Gegenbauer Polynomials}

In this section, we again deduce the closed form representation
for $A_n(t)$, this time utilizing an explicit representation of $w(z,t)$ 
in terms of Gegenbauer polynomials. Observe that this section is not
necessary for our development, but it shows some interesting connections
for the reader who is familiar with orthogonal polynomials and generating
functions.

Solving $w=K(z)=\frac{z}{(1-z)^2}$ for $z$ leads to the
representation
\[
K^{-1}(w)=
{\frac {1+2\,w-\sqrt {1+4\,w}}{2\,w}}
\]
for the inverse of the Koebe function. Therefore, substituting $e^{-t}K(z)$,
we obtain the representation
\be
w(z,t)=
{\frac {-1+(1+x)z-{z}^{2}+(1-z)\,\sqrt {1+{z}^{2}-2\,xz}}
{z\left (x-1\right )}}
\;,
\label{eq:wztfirst}
\ee
where we simplified the result changing variables according to
$e^{-t}=\frac{1-x}{2}$. 

Since $\sqrt {1+{z}^{2}-2\,xz}$ is the
generating function of the Gegenbauer polynomials $C_n^{(-1/2)}(x)$
(see e.\ g.\ \cite{AS}, (22.9.3)), (\ref{eq:wztfirst}) implies
for $n\geq 2$
\be
A_n(t)=\ed{x-1}\Big( C_{n+1}^{(-1/2)}(x)-C_{n}^{(-1/2)}(x)\Big)
\;.
\label{eq:Aktdifference}
\ee
On the other hand, it is well-known that $C_n^{(-1/2)}(x)$ has the 
(hypergeometric) representation
\be
C_n^{(-1/2)}(x)=2\,\sum_{j=0}^{n-1}\frac{(1-n)_j\,(n)_j}{j!\,(2)_j}
\lk\frac{1-x}{2}\rk^{j+1}
\label{eq:Gegenbauer-1/2}
\ee
when expanded at $x=1$, which can be obtained from \cite{AS}, (22.5.46),
as limiting case. Here, $(a)_j:=a(a+1)\cdots(a+j-1)$ as usual denotes the
{\sl Pochhammer symbol\/} (or {\sl shifted factorial\/}).

Therefore, we obtain for the difference (\ref{eq:Aktdifference})
\beao
A_n(t)&=&
\ed{x-1}\Big( C_{n+1}^{(-1/2)}(x)-C_{n}^{(-1/2)}(x)\Big)
\\&=&
-
\sum_{j=0}^{n}\frac{(-n)_j\,(n+1)_j}{j!\,(2)_j}
\lk\frac{1-x}{2}\rk^{j}
+
\sum_{j=0}^{n-1}\frac{(1-n)_j\,(n)_j}{j!\,(2)_j}
\lk\frac{1-x}{2}\rk^{j}
\\&=&
\sum_{j=1}^n \frac{(1-n)_{j-1}\,(n+1)_{j-1}}{2\,(j-1)!\,(2)_j}
\lk\frac{1-x}{2}\rk^{j}
=
n\,\frac{1-x}{2}\,\sum_{j=0}^{n-1}\frac{(1-n)_j\,(1+n)_j}{j!\,(3)_j}
\lk\frac{1-x}{2}\rk^{j}
\\&=&
n\,e^{-t}\,\sum_{j=0}^{n-1}\frac{(1-n)_j\,(1+n)_j}{j!\,(3)_j}\,e^{-jt}
=
\sum_{j=1}^n 2\,(-1)^{j+1}\ueber{n+j-1}{n-j}\frac{(2j-1)!}{(j-1)!(j+1)!}
\,e^{-jt}
\;.
\eeao
We note in passing
that the method presented in \cite{Koe92}--\cite{Koeortho} finds
the ordinary differential equation for $B_n(y)$ of 
Lemma~\ref{l:Ordinary differential equation for coefficient functions}, 
and furthermore a pure recurrence equation with respect to $n$,
automatically. Actually, this differential equation generated
by our {\sc Mathematica} \cite{Wol} implementation \cite{Koe93}
was an essential tool to discover the short proof of
Theorem~\ref{th:Coefficient representation of Lowner chain of Koebe function}.
Moreover, the same implementation discovers the power series
representation (\ref{eq:Gegenbauer-1/2}) automatically.

\section{The de Branges and Weinstein functions}
\label{sec:The de Branges and Weinstein functions}

In \cite{Bra} de Branges showed that the Milin conjecture is valid
if for all $n\geq 2$ the {\sl de Branges functions} 
$\tau_{k}^n:\R^+\rightarrow\R
\;(k=1,\ldots,n+1)$ defined by the system of differential equations
\bea
\tau_{k+1}^n(t)-\tau_{k}^n(t)&=&
\displaystyle{\frac{\textstyle{\dot{\tau}_{k}^n(t)}}{\textstyle{k}}}
+\displaystyle{\frac{\textstyle{\dot{\tau}_{k+1}^n(t)}}{\textstyle{k+1}}}
\;\;\;\;\;\;\kn
\label{eq:dB2}
\\[4mm]
\tau_{n+1}^n&\equiv& 0
\label{eq:dB2a}
\eea
with the initial values
\be
\tau_{k}^n(0)=n+1-k
\label{eq:dB1}
\ee
have the properties
\be
\lim\limits_{t\rightarrow\infty}\tau_{k}^n(t)=0\;,
\label{eq:dB3}
\ee
and
\be
\dot{\tau_{k}^n}(t)\leq 0\;\;\;\;\;(t\in\R^+)
\;.
\label{eq:dB4}
\ee
The relation (\ref{eq:dB3}) is easily checked using standard methods
for ordinary differential equations, whereas (\ref{eq:dB4}) is a deep result.

L.\ de Branges gave an explicit representation of the function system
$\tau_k^n(t)$ (\cite{Bra}, \cite{Henrici}, \cite{Schmersau})
(that we don't use, though, see 
\S~\ref{sec:Closed Form Representation of Weinstein functions},
however), with which the proof of the 
de Branges theorem was completed as soon as de Branges realized that
(\ref{eq:dB4}) was a theorem previously proved by Askey and Gasper \cite{AG}.

Note that the derivatives $\dot\tau_k^n(t)$ are characterized by
the same system of differential equations (\ref{eq:dB2}), the equation
\be
\dot\tau_n^n(t)=-n\,e^{-nt}
\label{eq:dottaunn}
\ee
and the initial values
\be
\dot\tau_k^n(0)=\funkdeff{-k}{n-k\;\;\;\mbox{even}}{0}{n-k\;\;\;\mbox{odd}}
\label{eq:taupunktIV}
\ee
as replacements for (\ref{eq:dB2a}) and
(\ref{eq:dB1}) (see e.\ g.\ \cite{FP}, p.\ 685).

On the other hand, Weinstein \cite{Weinstein} 
uses the L\"owner chain (\ref{eq:Koebe}),
and shows the validity of Milin's conjecture if for all $n\geq 2$ 
the {\sl Weinstein functions} $\Lambda_k^n:\R^+\rightarrow\R
\;(k=1,\ldots,n+1)$ defined by
\be
\frac{e^{t}w(z,t)^{k+1}}{1-w^2(z,t)}=:
\sum\limits_{n=k}^{\infty}\Lambda_k^n(t)z^{n+1}
=W_k(z,t)
\label{eq:W(z,t)}
\ee
satisfy the relations
\begin{equation}
\Lambda_k^n(t)\geq 0\;\;\;\;\;\;(t\in\R^+,\;\;\;k,n\in\N)
\;.
\label{eq:pos}
\end{equation}
Weinstein did not identify the functions $\Lambda_k^n(t)$, but was
able to prove (\ref{eq:pos}) without an explicit representation.

Independently, both Todorov \cite{Todorov} and Wilf \cite{Wilf} 
proved---using the explicit representation of the de Branges 
functions---the following
\bT
\label{th:Connection between de Branges and Weinstein functions}
{\bf (Connection between de Branges and Weinstein functions)}
{\rm
For all\linebreak
$n\in\N$, $k=1,\ldots,n$, one has the identity
\be
\dot{\tau_{k}^n}(t)=-k\Lambda_k^n(t)
\;,
\label{eq:whatweliketohave}
\ee
i.\ e.\ the de Branges and the Weinstein functions essentially are the same,
and the main inequalities (\ref{eq:dB4}) and (\ref{eq:pos}) are identical.
\eop
}
\eT
In this section, we give a very elementary proof of this result, which 
in view of
(\ref{eq:W(z,t)}) can be looked at as a property of the Koebe function.

Firstly, we realize that (again, we omit the arguments of $w(z,t)$)
\be
W_{k+1}(z,t)=\frac{e^{t}w^{k+2}}{1-w^{2}}=w\,W_{k}(z,t)
\;,
\label{eq:WkintermsofWk-1}
\ee
and that further
\[
W_{k}(z,t)=\frac{e^{t}w^{k+1}}{1-w^{2}}=\frac{e^{t}w}{(1-w)^2}
\frac{1-w}{1+w}w^{k}=
K(z)\frac{1-w}{1+w}w^{k}
\;,
\]
so that with (\ref{eq:w(z,t)LoewnerDE}) in particular
\be
W_{1}(z,t)=
K(z)\frac{1-w}{1+w}w
=
-K(z)\,\dot w
\;.
\label{eq:WeinsteinAnfang}
\ee
Moreover, we get the relation
\[
W_{k}(z,t)+W_{k+1}(z,t)=(1+w)W_{k}(z,t)=K(z)(1-w)w^{k}=
K(z)w^k-K(z)w^{k+1}
\;.
\]
Taking derivative with respect to $t$, this identity implies
\bea
\dot W_{k}(z,t)+\dot W_{k+1}(z,t)
&=&
K(z)\cdot k\,w^{k-1}\,\dot w-K(z)\cdot (k+1)\,w^{k}\,\dot w
\label{eq:drefor wk}
\\&=&
-k\cdot K(z)\,\frac{1-w}{1+w}\,w^{k}+(k+1)\cdot K(z)\,\frac{1-w}{1+w}\,w^{k+1}
\nonumber
\\&=&
(k+1)\,W_{k+1}(z,t)-k\,W_{k}(z,t)
\nonumber
\eea
where again,
we utilized the L\"owner differential equation (\ref{eq:w(z,t)LoewnerDE})
for $w(z,t)$.

Equating coefficients it follows that the same system of 
differential equations is valid for $\Lambda_k^n(t)$, and therefore
for $y_k^n(t):=-k\,\Lambda_k^n(t)$ we get the differential equations
system
\[
y_{k+1}^n(t)-y_k^n(t)=
\frac{\dot y_k^n(t)}{k}+\frac{\dot y_{k+1}^n(t)}{k+1}
\]
of de Branges (\ref{eq:dB2}).  From
\[
\Lambda_n^n(t)=
\lim_{z\pf 0}\frac{W_n(z,t)}{z^n}
=
\lim_{z\pf 0}\frac{e^t}{1-w^2(z,t)}\lk\frac{w(z,t)}{z}\rk^{n+1}
=
e^t\lk e^{-(n+1)t}\rk=
e^{-nt}
\]
which follows from (\ref{eq:W(z,t)}), we realize that
\[
y_n^n(t)=-n\,e^{-nt}
\]
so that (\ref{eq:dottaunn}) is satisfied.

To show (\ref{eq:whatweliketohave}), it therefore remains to prove
(\ref{eq:taupunktIV}) which can be read off from 
\[
W_k(z,0)=\sum_{n=k}^\infty \Lambda_k^n(0)z^{n+1}=\frac{z^k}{1-z^2}
=\sum_{j=0}^\infty z^{2j+k}
\;.
\]
This finishes the proof of 
Theorem~\ref{th:Connection between de Branges and Weinstein functions}.

\section{Closed Form Representation of 
\mbox{Weinstein functions}}
\label{sec:Closed Form Representation of Weinstein functions}

In this section, we show how---in a similar manner as we derived the closed
form representation of the coefficients of Koebe's L\"owner chain $w(z,t)$ 
in \S~\ref{sec:The Lowner Chain of the Koebe Function}---the closed form 
representation of $\dot{\tau_{k}^n}(t)=-k\Lambda_k^n(t)$ that was
given by de Branges, can be deduced in an elementary way,
only utilizing the properties of $w(z,t)$ 
that we developed in \S~\ref{sec:The Lowner Chain of the Koebe Function}.
In particular, the known proofs of the Bieberbach and Milin conjectures  
may be regarded as a consequence
of the L\"owner differential equation, plus properties of the Koebe function.

Since $W_1(z,t)$ is given by (\ref{eq:WeinsteinAnfang}) in terms of
$w(z,t)$, and $W_k(z,t)$
satisfies the recurrence (\ref{eq:WkintermsofWk-1}),
from the representation (\ref{eq:closedform:w(z,t)}) of $w(z,t)$ we deduce
by induction that the coefficients $\Lambda_k^n(t)$ of $W_k(z,t)$
have a representation
\be
\Lambda_k^n(t)=
\sum_{j=k}^n a_j^{(n,k)}\,e^{-jt}
\quad\quad(n\geq k)
\;.
\label{eq:ansatzlambdakn}
\ee
Substituting $\Lambda_k^n(t)$ according to
(\ref{eq:W(z,t)}) in (\ref{eq:drefor wk}), 
and equating coefficients, we obtain for $n\geq k+1\geq 2$ 
\be
\dot \Lambda_k^n(t)+\dot\Lambda_{k+1}^n(t)
=
(k+1)\Lambda_{k+1}^n(t)-k\Lambda_k^n(t)
\;.
\label{eq:lambdapunktrelation}
\ee
If we substitute now (\ref{eq:ansatzlambdakn}) 
in (\ref{eq:lambdapunktrelation}), and equate coefficients, again, then
we get the simple recurrence equation $n\geq j\geq k\geq 2$
\[
a_j^{(n,k)}=-\frac{j-k+1}{j+k}a_j^{(n,k-1)}
\]
for the coefficients $a_j^{(n,k)}$ which (by telescoping) generates
\be
a_j^{(n,k)}=(-1)^{k-1}\,\frac{(j-1)!\,(j+1)!}{(j-k)!\,(j+k)!}\,
a_j^{(n,1)}
\quad\quad
(n\geq j\geq k\geq 2)
\;.
\label{eq:ajnk}
\ee
Therefore, to get a closed form representation of $a_j^{(n,k)}$,
we need only one for $a_j^{(n,1)}$, and we are done. To obtain this result,
we observe that
\[
\Lambda_1^n(t)=-\sum_{l=1}^n (n+1-l)\,\dot A_l(t)
\]
following from (\ref{eq:A(k,t)}) and (\ref{eq:WeinsteinAnfang}).
Using the definition
$a_j^{(n,k)}:= 0$ for $j<k$, and the representation (\ref{eq:Akfineal})
of $A_l(t)$ that we gave in
Theorem \ref{th:Coefficient representation of Lowner chain of Koebe function},
we deduce in a straightforward manner
\beao
\sum_{j=1}^n a_j^{(n,1)}\,e^{-jt}
&=&
\Lambda_1^n(t) = -\sum_{l=1}^n (n+1-l)\,\dot A_l(t)
\\&=&
\sum_{l=1}^n (n+1-l)\sum_{j=1}^n 
2j(-1)^{j+1}\ueber{l+j-1}{l-j}\frac{(2j-1)!}{(j-1)!(j+1)!} \,e^{-jt}
\\&=&
\sum_{j=1}^n 
2j(-1)^{j+1}\frac{(2j-1)!}{(j-1)!(j+1)!}
\lk
\sum_{l=1}^n
(n+1-l)\ueber{l+j-1}{l-j}
\rk
e^{-jt}
\eeao
where we changed the order of summation. Equating coefficients, we 
therefore see that
\[
a_j^{(n,1)}=
2j(-1)^{j+1}\frac{(2j-1)!}{(j-1)!(j+1)!}
\lk
\sum_{l=j}^n
(n+1-l)\ueber{l+j-1}{l-j}
\rk
\;.
\]
Since for $b_l:=(n+1-l)\ueber{l+j-1}{l-j}$, one has
\[
b_l=s_l-s_{l-1}
\quad\quad
\mbox{with}
\quad\quad
s_l:=\frac{(j+l)(n+1+j+2jn-2jl)}{2j(2j+1)}\ueber{l+j-1}{l-j}
\;,
\]
(i.\ e., $s_l$ is an antidifference of $b_l$ which is found by 
Gosper's algorithm, see \cite{Gosper}, \cite{Koepf95})
which can easily be checked, and since $s_{j-1}\equiv 0$ it turns out that
\beao
\sum_{l=j}^n b_l&=&
\sum_{l=j}^n (n+1-l)\ueber{l+j-1}{l-j}
\\&=&
\sum_{l=j}^n (s_l-s_{l-1})
=
s_n-s_{j-1}
=
\frac{(j+n)(n+1+j)}{2j(2j+1)}\ueber{n+j-1}{n-j}
\;.
\eeao
Therefore, using (\ref{eq:ajnk}), we finally have 
\bea
a_j^{(n,k)}&=&(-1)^{k-1}\,\frac{(j-1)!\,(j+1)!}{(j-k)!\,(j+k)!}\,
a_j^{(n,1)}
\nonumber
\\&=&
(-1)^{k-1}\,\frac{(j-1)!\,(j+1)!}{(j-k)!\,(j+k)!}\,
2j(-1)^{j+1}\frac{(2j-1)!}{(j-1)!(j+1)!}
\lk
\sum_{l=j}^n
(n+1-l)\ueber{l+j-1}{l-j}
\rk
\nonumber
\\&=&
(-1)^{k+j}\,\frac{(2j-1)!}{(j-k)!\,(j+k)!}\,
\frac{(j+n)(n+1+j)}{(2j+1)}\ueber{n+j-1}{n-j}
\nonumber
\\&=&
(-1)^{k+j}\,\ueber{2j}{j-k}\ueber{n+j+1}{n-j}
\;,
\label{eq:finalajnk}
\eea
and hence 
\be
\Lambda_k^n(t)=
\sum_{j=k}^n a_j^{(n,k)}\,e^{-jt}
=
\sum_{j=k}^n
(-1)^{k+j}\,\ueber{2j}{j-k}\ueber{n+j+1}{n-j}
\,e^{-jt}
\;. 
\label{eq:finalLambda}
\ee
This, by (\ref{eq:whatweliketohave}), gives de Branges'
closed representation for $\dot{\tau_{k}^n}(t)$.

Note that in our presentation no knowledge about hypergeometric
functions is needed. On the other hand, from representation
(\ref{eq:finalajnk}) one can read off
\[
\frac{a_{j+1}^{(n,k)}}{a_j^{(n,k)}}
=
{\frac {\left (j+1/2\right )\left (j+n+2\right )\left (j-n\right )}{
\left (j+3/2\right )\left (j+k+1\right )\left (j-k+1\right )}}
\;,
\]
and since $k$ is an integer, we may substitute $j\mapsto j+k$
(i.\ e.\ shift the summation variable) which leads to
\[
\frac{a_{j+k+1}^{(n,k)}}{a_{j+k}^{(n,k)}}
=
{\frac {\left (j+k+1/2\right )\left (j+n+k+2\right )\left (j-n+k
\right )}{\left (j+k+3/2\right )\left (j+2\,k+1\right )\left (j+1
\right )}}
\]
and to the initial value
\[
e^{-kt}\,a_k^{(n,k)}=
e^{-kt}\ueber{n+k+1}{n-k}
\]
so that (\ref{eq:finalLambda}) reads
\[
\Lambda_k^n(t)=
e^{-kt}\ueber{n+k+1}{n-k}
\;
_3 F_2\left.
\!\!
\left(
\!\!\!\!
\begin{array}{c}
\multicolumn{1}{c}{\begin{array}{c}
k+1/2\;, n+k+2\;, -n+k
\end{array}}\\[1mm]
\multicolumn{1}{c}{\begin{array}{c}
k+3/2\;,2\,k+1
            \end{array}}\end{array}
\!\!\!\!
\right| e^{-t}\right)
\;.
\]
Note, that similarly from (\ref{eq:Akfineal}) 
one gets the hypergeometric representations
\[
A_n(t)=
n\,e^{-t}
\;
_2 F_1\left.
\!\!
\left(
\!\!\!\!
\begin{array}{c}
\multicolumn{1}{c}{\begin{array}{c}
1-n\;, n+1
\end{array}}\\[1mm]
\multicolumn{1}{c}{\begin{array}{c}
3
            \end{array}}\end{array}
\!\!\!\!
\right| e^{-t}\right)
\;,
\]
and for the Gegenbauer polyomials $C_n^{(-1/2)}(x)$ by (\ref{eq:Gegenbauer-1/2})
\[
C_n^{(-1/2)}(x)=
(1-x)
\;
_2 F_1\left.
\!\!
\left(
\!\!\!\!
\begin{array}{c}
\multicolumn{1}{c}{\begin{array}{c}
1-n\;, n
\end{array}}\\[1mm]
\multicolumn{1}{c}{\begin{array}{c}
2
            \end{array}}\end{array}
\!\!\!\!
\right| \frac{1-x}{2}\right)
\;.
\]

\section{Generating Function of the de Branges Functions}

In this final section, we give a very simple representation of the
generating function $B_k(z,t)$ of the de Branges functions
\[
B_k(z,t)=\sum_{n=k}^\infty \tau^n_k(t)\,z^{n+1}
\]
from which one can directly deduce de Branges' main inequality
$\tau^n_k(t)\geq 0$ (see \cite{Bra}, \cite{FP})
without utilizing the inequality for the
derivatives $\dot{\tau^n_k}(t)\leq 0$.

Whereas de Branges considered the Milin conjecture for {\sl fixed} $n\in\N$,
and therefore introduced $\tau_k^n(t)\;
(k=1,\ldots,n+1)$, we take a fixed $k\in\N$ and the
generating function of $\tau^n_k(t)$ with respect to $n$, 
hence all $n\geq k$ are considered {\sl at the same time}.
\bT
{\rm
The generating function of the de Branges functions has the
representation
\be
B_k(z,t)=\sum_{n=k}^\infty \tau^n_k(t)\,z^{n+1}=
K(z)\,w(z,t)^k=
K(z)\,\lk\frac{4e^{-t} z}{\left(1-z+\sqrt{1-2xz+z^2}\right)^2}\rk^k
,
\label{eq:second representation}
\ee
$(x=1-2e^{-t})$. Moreover one has the hypergeometric
representation
\bea
B_k(z,t)
&=&
K(z)^{k+1}\,e^{-kt}\;\hypergeom{2}{1}{k,k+1/2}{2k+1}{-4K(z)e^{-t}}
\nonumber
\\&=&
\sum_{j=k}^\infty 
(-1)^{j+k}\,\frac{2k}{j+k}\,\ueber{2j-1}{j-k}\,K(z)^{j+1}\,e^{-jt}
\label{eq:explicitgendeB}
\;,
\eea
being a Taylor series representation with respect to $y=e^{-t}=\frac{1-x}{2}$.
}
\eT
\pr
Define $B_k(z,t)$ by
\be
B_k(z,t):=K(z)\,w(z,t)^k
\;;
\label{eq:Bkdef}
\ee
then 
\[
\dot{B_{k}}(z,t)=K(z)\,k\,w(z,t)^{k-1}\,\dot{w}(z,t)
\;.
\]
Using (\ref{eq:w(z,t)LoewnerDE}), we get therefore
\bea
\frac{\dot{B}_{k+1}(z,t)}{k+1}+\frac{\dot{B_{k}}(z,t)}{k}
\nonumber
&=&
K(z)\,w(z,t)^{k-1}\,\dot{w}(z,t)\,(1+w(z,t))
\\&=&
-K(z)\,w(z,t)^{k}(1-w(z,t))=K(z)\,w(z,t)^{k+1}-K(z)\,w(z,t)^{k}
\nonumber
\\&=&
B_{k+1}(z,t)-B_k(z,t)
\;.
\label{eq:recBk}
\eea
By the definition (\ref{eq:Bkdef}) of $B_k(z,t)$ its Taylor series for
$z=0$ starts with a $z^{k+1}$ term, hence we may write
\be
B_k(z,t)=\sum\limits_{n=k}^\infty \tau^n_k(t)\,z^{n+1}
\label{eq:Bksum}
\ee
and we can assume $\tau^n_{n+1}(t)\equiv 0$, hence (\ref{eq:dB2a}).
Substituting (\ref{eq:Bksum}) in (\ref{eq:recBk})
yields furthermore (\ref{eq:dB2}) by equating coefficients of $z^{n+1}$.
The $(n+1)^{st}$ Taylor coefficient of
\[
B_k(z,0)=K(z)\,w(z,0)^k=\frac{z^{k+1}}{(1-z)^2}
\]
equals $n+1-k$,
hence the initial values (\ref{eq:dB1}) are satisfied, and therefore
$\tau^n_k(t)$ form the de Branges functions.

Starting with (\ref{eq:wztfirst}),
a calculation shows that $w(z,t)$ has the explicit representation
\[
w(z,t)=\frac{4e^{-t}z}{\lk 1-z+\sqrt{1-2xz+z^2}\rk^2}
\;,
\]
hence the right hand representation of (\ref{eq:second representation}) 
follows.

In a similar manner as we derived the closed
form representation of the coefficients of Koebe's L\"owner chain $w(z,t)$
in \S~\ref{sec:The Lowner Chain of the Koebe Function} (or
by the method presented in \cite{Koe92}--\cite{Koeortho}), one deduces
(\ref{eq:explicitgendeB}).
\eop
To deduce the inequalities $\tau^n_k(t)\geq 0$ as announced, we remark that
the Jacobi polynomials $P_n^{(\al,\beta)}(x)$ have the generating function
\[
\sum_{n=0}^\infty P_n^{(\al,\beta)}(x)\,z^n=
\frac{2^{\al+\beta}}{\sqrt{1\!-\!2xz\!+\!z^2}}\,
\frac{1}{\left(1-z+\sqrt{1\!-\!2xz\!+\!z^2}\right)^{\alpha}}\,
\frac{1}{\left(1+z+\sqrt{1\!-\!2xz\!+\!z^2}\right)^{\beta}}
\]
(see e.g.\ \cite{AS}, (22.9.1)), hence
\beao
B_k(z,t)
&=&
K(z)\,\lk\frac{4e^{-t} z}{\left(1-z+\sqrt{1-2xz+z^2}\right)^2}\rk^k
\\&=&
z^{k+1}\,e^{-kt} \cdot\ed{1-z}\,
\frac{2^{2k}}{\sqrt{1\!-\!2xz\!+\!z^2}}\,
\frac{1}{\left(1-z+\sqrt{1\!-\!2xz\!+\!z^2}\right)^{2k}}
\cdot \frac{\sqrt{1\!-\!2xz\!+\!z^2}}{1-z}
\\&=&
z^{k+1}\,e^{-kt} \cdot
\sum_{n=0}^\infty\sum_{j=0}^n P_j^{(2k,0)}(x)\,z^n
\cdot
\sum_{n=0}^\infty\sum_{j=0}^n C_j^{(-1/2)}(x)\,z^n
\;,
\eeao
and the result follows from the positivity of the Jacobi polynomial sums
\be
\sum_{j=0}^{n} P_j^{(2k,0)}(x)\geq 0
\label{eq:AskeyGasper}
\ee
(see \cite{AG}, Theorem 3) and the positivity of the Taylor coefficients 
(with respect to $z$) of the function 
\[
\frac{\sqrt{1\!-\!2xz\!+\!z^2}}{1-z}
\]
(see \cite{AG}, Theorem D). Hence, again, as in
de Branges' original proof, the Askey-Gasper result (\ref{eq:AskeyGasper})
does the main job.


\end{document}